\documentclass[11pt]{article}
\usepackage[utf8]{inputenc}
\usepackage[utf8]{inputenc}
\usepackage{hyperref}
\usepackage{amsmath,amssymb,amsthm}
\usepackage{amsfonts,tikz}
\usepackage[margin=1in]{geometry}
\usepackage{mathtools}
\usetikzlibrary{positioning,automata}
\tikzset{every loop/.style={min distance=10 mm, in=60, out=120, looseness=10}}
\newcommand{\bk}{\backslash}
\newcommand{\lm}{\lambda}

\newcommand{\mg}{\mathbb{G}}

\newcommand{\tm}{\widetilde{m}}

\newtheorem{lemma}[]{Lemma}
\newtheorem{cor}[lemma]{Corollary}
\newtheorem{theorem}[lemma]{Theorem}

\newtheorem{definition}[lemma]{Definition}

\title{Modular Relations of the Tutte Symmetric Function}
\author{Logan Crew, Sophie Spirkl\footnote{Department of Combinatorics \& Optimization, University of Waterloo, Waterloo, ON, N2L 3E9.\newline  Emails:  lcrew@uwaterloo.ca, sspirkl@uwaterloo.ca. \newline
We acknowledge the support of the Natural Sciences and Engineering Research Council of Canada (NSERC), [funding
reference number RGPIN-2020-03912]. \newline Cette recherche a \'et\'e financ\'ee par le Conseil de recherches en sciences naturelles et en g\'enie du Canada (CRSNG),
[num\'ero de r\'ef\'erence RGPIN-2020-03912]. }}
\date{\today}

\begin{document}

\maketitle

\begin{abstract}

For a graph $G$, its Tutte symmetric function $XB_G$ generalizes both the Tutte polynomial $T_G$ and the chromatic symmetric function $X_G$. We may also consider $XB$ as a map from the $t$-extended Hopf algebra $\mathbb{G}[t]$ of labelled graphs to symmetric functions.

We show that the kernel of $XB$ is generated by vertex-relabellings and a finite set of modular relations, in the same style as a recent analogous result by Penaguiao on the chromatic symmetric function $X$. In particular, we find one such relation that generalizes the well-known triangular modular relation of Orellana and Scott, and build upon this to give a modular relation of the Tutte symmetric function for any two-edge-connected graph that generalizes the $n$-cycle relation of Dahlberg and van Willigenburg. Additionally, we give a structural characterization of all local modular relations of the chromatic and Tutte symmetric functions, and prove that there is no single local modification that preserves either function on simple graphs.
    
\end{abstract}

\section{Introduction}

The \emph{Tutte symmetric function}, introduced by Stanley \cite{stanley2}, is a simultaneous extension of the Tutte polynomial and the chromatic symmetric function of a graph. Both of these latter functions have been and continue to be very well-studied. The Tutte polynomial is a universal deletion-contraction polynomial for graphs, containing information on many important functions such as the chromatic and flow polynomials, and very recently an entire handbook has been published cataloguing and chronicling its relevant applications to graphs, and more generally to matroids \cite{jo}. The chromatic symmetric function extends the chromatic polynomial by providing further information on graphs, such as counting acyclic orientations by number of sinks \cite{stanley}, while also providing deep connections to algebraic geometry \cite{unit, categ, wachs}. Current research on the chromatic symmetric function focuses on the Stanley-Stembridge conjecture that unit interval graphs are $e$-positive \cite{huh, epos, dahl, foley}, and the conjecture that the function distinguishes nonisomorphic trees \cite{trees, heil}. 

In contrast, the Tutte symmetric function has received very little attention since its creation. Although it has been shown to be equivalent to other graph functions such as the $W$-polynomial \cite{noble}, the polychromate \cite{merino, sarm}, and the $(r,q)$-chromatic function \cite{potts}, little is known about its expansions in symmetric function bases or what structural information it encodes about a graph beyond what the chromatic symmetric function already does. Recently the authors, Aliste-Prieto, and Zamora considered a vertex-weighted version \cite{tutteCrew}, and used this to find a deletion-contraction relation and a spanning tree formula that are reminiscent of similar well-known results for the Tutte polynomial. In this same work, the authors also explored how to construct pairs of graphs with equal Tutte symmetric function.

In this paper, we build on \cite{tutteCrew}. After giving background on symmetric function theory and graph colorings in Section 2, we proceed in Section 3 to consider $X$ and $XB$ as maps from formal linear combinations of vertex-labelled graphs to symmetric functions, and give a structural characterization of all combinations in $Ker(XB)$ that are universal in the sense that they remain in $Ker(XB)$ even upon adding arbitrary vertices and edges uniformly to all graphs in the combination. In other words, we characterize all ways to locally modify subgraphs of any larger graph in a way that preserves $XB$. We also provide such a characterization for $Ker(X)$, which allows us to prove that there is no single local modification that will transform a graph into another distinct graph with the same chromatic symmetric function (much less Tutte symmetric function).

Finally, in Section 4, we determine $Ker(XB)$ as the span of an explicit set of generators; equivalently, we give a finite number of local modifications such that any two graphs with the same Tutte symmetric function are related by a finite number of these modifications. Notably, one of these relations is a generalization of the triangular relation of Orellana and Scott \cite{ore}, which may be used to derive a simple relation for $X$ or $XB$ that is applicable to a broader range of cases. Additionally, we give a related modular relation that applies to any two-edge-connected graph, yielding a result that also generalizes the $n$-cycle modular relation for the chromatic symmetric function given by Dahlberg and van Willigenburg \cite{dahl}.
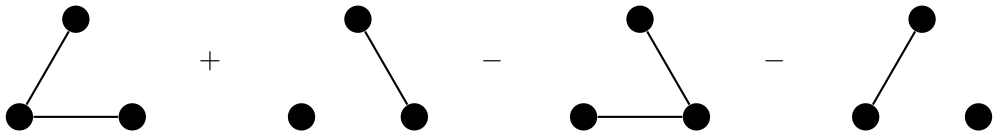
\begin{figure}[hbt]
\begin{center}
  \begin{tikzpicture}[scale=1.5]
  \node[fill=black, circle] at (0, 0)(1){};
    \node[fill=black, circle] at (1, 0)(2){};
    \node[fill=black, circle] at (0.5, 0.866)(3){};

    \draw[black, thick] (3) -- (1);
    \draw[black, thick] (2) -- (1);

    \draw (1.5, 0.5) coordinate (PL) node[right] { $\bf{+}$ };
    
    \node[fill=black, circle] at (2.5, 0)(4){};
    \node[fill=black, circle] at (3.5, 0)(5){};
    \node[fill=black, circle] at (3, 0.866)(6){};
    
    \draw[black, thick] (5) -- (6);
    
    \draw (4, 0.5) coordinate (MI) node[right] { $\bf{-}$ };
    
    \node[fill=black, circle] at (5, 0)(7){};
    \node[fill=black, circle] at (6, 0)(8){};
    \node[fill=black, circle] at (5.5, 0.866)(9){};
    
    \draw[black, thick] (8) -- (7);
    \draw[black, thick] (9) -- (8);
    
    \draw (6.5, 0.5) coordinate (M2) node[right] { $\bf{-}$ };
    
    \node[fill=black, circle] at (7.5, 0)(A){};
    \node[fill=black, circle] at (8.5, 0)(B){};
    \node[fill=black, circle] at (8, 0.866)(C){};
    
    \draw[black, thick] (A) -- (C);
  
  \end{tikzpicture}
\end{center}
\label{fig:newfig}
\caption{A relation in the kernel of $XB$ also valid for $X$}
\end{figure}

\section{Background}

\subsection{Fundamentals of Partitions and Symmetric Functions}

 A \emph{set partition} of a set $S$ is a collection of nonempty, pairwise nonintersecting \emph{blocks} $B_1, \dots, B_k$ satisfying $B_1 \cup \dots \cup B_k = S$. We will specify that a union of blocks is a set partition by writing $\sqcup$ for disjoint union, using the notation $B_1 \sqcup \dots \sqcup B_k$. We will frequently consider partitions of the set $\{1,2,\dots,n\} = [n]$. When $n$ is clear, for shorthand, we will write $p(i_{11} \dots i_{1k_1}, i_{21} \dots i_{2k_2}, $ $\dots, i_{m1} \dots i_{mk_m})$ to denote the partition of $[n]$ whose blocks are all singletons except for $\{i_{11}, \dots i_{1k_1}\},$ $ \{i_{21} \dots i_{2k_2}\}, \dots, $ $\{i_{m1}, \dots, i_{mk_m}\}$. For example, when $n = 6$, we write $p(13,26)$ to mean the set partition $\{1,3\} \sqcup \{2,6\} \sqcup \{4\} \sqcup \{5\}$.
 
 An \emph{integer partition} is a tuple $\lm = (\lm_1,\dots,\lm_k)$ of positive integers such that $\lm_1 \geq \dots \geq \lm_k$.  The integers $\lm_i$ are the \emph{parts} of $\lm$.  If $\sum_{i=1}^k \lm_i = n$, we say that $\lm$ is a partition of $n$. The number of parts equal to $i$ in $\lm$ is given by $r_i(\lm)$.
 
 We may use simply \emph{partition} to refer to either a set or integer partition. We write $\pi \vdash [n]$ or $\lm \vdash n$ to mean respectively that $\pi$ is a partition of $[n]$, and $\lm$ is a partition of $n$, and we write $|\lm| = |\pi| = n$. The number of blocks or parts is the \emph{length} of a partition, and is denoted by $l(\pi)$ or $l(\lm)$. When $\pi$ is a set partition, we will write $\lm(\pi)$ to mean the integer partition whose parts are the sizes of the blocks of $\pi$.

  A function $f(x_1,x_2,\dots) \in \mathbb{C}[[x_1,x_2,\dots]]$ is \emph{symmetric}\footnote{The choice of coefficient ring is irrelevant for the work in this paper so long as it is a field of characteristic $0$.} if $f(x_1,x_2,\dots) = f(x_{\sigma(1)},x_{\sigma(2)},\dots)$ for every permutation $\sigma$ of the positive integers $\mathbb{N}$.  The \emph{algebra of symmetric functions} $\Lambda$ is the subalgebra of $\mathbb{C}[[x_1,x_2,\dots]]$ consisting of those symmetric functions $f$ that are of bounded degree (that is, there exists a positive integer $n$ such that every monomial of $f$ has degree $\leq n$).  Furthermore, $\Lambda$ is a graded algebra, with natural grading
  $$
  \Lambda = \bigoplus_{d=0}^{\infty} \Lambda^d
  $$
  where $\Lambda^d$ consists of symmetric functions that are homogeneous of degree $d$ \cite{mac,stanleybook}.

  Each $\Lambda^d$ is a finite-dimensional vector space over $\mathbb{C}$, with dimension equal to the number of partitions of $d$ (and thus, $\Lambda$ is an infinite-dimensional vector space over $\mathbb{C}$).  Some commonly-used bases of $\Lambda$ that are indexed by partitions $\lm = (\lm_1,\dots,\lm_k)$ include:
\begin{itemize}
  \item The monomial symmetric functions $m_{\lm}$, defined as the sum of all distinct monomials of the form $x_{i_1}^{\lm_1} \dots x_{i_k}^{\lm_k}$ with distinct indices $i_1, \dots, i_k$.

  \item The power-sum symmetric functions, defined by the equations
  $$
  p_n = \sum_{k=1}^{\infty} x_k^n, \hspace{0.3cm} p_{\lm} = p_{\lm_1}p_{\lm_2} \dots p_{\lm_k}.
  $$
  \item The elementary symmetric functions, defined by the equations
  $$
  e_n = \sum_{i_1 < \dots < i_n} x_{i_1} \dots x_{i_n}, \hspace{0.3cm} e_{\lm} = e_{\lm_1}e_{\lm_2} \dots e_{\lm_k}.
  $$
\end{itemize}

  We also make use of the \emph{augmented monomial symmetric functions}, defined by 
  $$
  \tm_{\lm} = \left(\prod_{i=1}^{\infty} r_i(\lm)!\right)m_{\lm}.
  $$
  
  Given a symmetric function $f$ and a basis $b$ of $\Lambda$, we say that $f$ is \emph{$b$-positive} if when we write $f$ in the basis $b$, all coefficients are nonnegative.
  
  \subsection{Fundamentals of Graphs and Colorings}

  A \emph{graph} $G = (V,E)$ consists of a \emph{vertex set} $V$ and an \emph{edge multiset} $E$ where the elements of $E$ are (unordered) pairs of (not necessarily distinct) elements of $V$. An edge $e \in E$ that contains the same vertex twice is called a \emph{loop}.  If there are two or more edges that each contain the same two vertices, they are called \emph{multi-edges}.  A \emph{simple graph} is a graph $G = (V,E)$ in which $E$ does not contain loops or multi-edges (thus, $E \subseteq \binom{V}{2}$).  If $\{v_1,v_2\}$ is an edge (or nonedge), we will write it as $v_1v_2 = v_2v_1$.  The vertices $v_1$ and $v_2$ are the \emph{endpoints} of the edge $v_1v_2$. We will use $V(G)$ and $E(G)$ to denote the vertex set and edge multiset of a graph $G$, respectively.
  
  Two graphs $G$ and $H$ are said to be \emph{isomorphic} if there exists a bijective map $f: V(G) \rightarrow V(H)$ such that for all $v_1,v_2 \in V(G)$ (not necessarily distinct), the number of edges $v_1v_2$ in $E(G)$ is the same as the number of edges $f(v_1)f(v_2)$ in $E(H)$.

  The \emph{complement} of a simple graph $G = (V,E)$ is denoted $\overline{G}$, and is defined as $\overline{G} = (V, \binom{V}{2} \bk E)$, so in $\overline{G}$ every edge of $G$ is replaced by a nonedge, and every nonedge is replaced by an edge.
  
  For a vertex $v \in V(G)$, its \emph{degree} $d(v)$ is the number of times $v$ occurs as an endpoint of an edge in $E(G)$ (thus a loop at $v$ adds $2$ to the degree).

  A \emph{subgraph} of a graph $G$ is a graph $G' = (V',E')$ where $V' \subseteq V$ and $E' \subseteq E|_{V'}$, where $E|_{V'}$ is the set of edges with both endpoints in $V'$.  An \emph{induced subgraph} of $G$ is a graph $G' = (V',E|_{V'})$ with $V' \subseteq V$.  The induced subgraph of $G$ using vertex set $V'$ will be denoted $G|_{V'}$.  A \emph{stable set} of $G$ is a subset $V' \subseteq V$ such that $E|_{V'} = \emptyset$.  A \emph{clique} of $G$ is a subset $V' \subseteq V$ such that for every pair of distinct vertices $v_1$ and $v_2$ of $V'$, $v_1v_2 \in E(G)$.

  A \emph{path} in a graph $G$ is a nonempty sequence of edges $v_1v_2$, $v_2v_3$, \dots, $v_{k-1}v_k$ such that $v_i \neq v_j$ for all $i \neq j$.  The vertices $v_1$ and $v_k$ are the \emph{endpoints} of the path. A \emph{cycle} in a graph is a nonempty sequence of distinct edges $v_1v_2$, $v_2v_3$, \dots, $v_kv_1$ such that $v_i \neq v_j$ for all $i \neq j$. Note that in a simple graph every cycle must have at least $3$ edges, although in a nonsimple graph there may be cycles of size $1$ (a loop) or $2$ (multi-edges).

  A graph $G$ is \emph{connected} if for every pair of vertices $v_1$ and $v_2$ of $G$ there is a path in $G$ with $v_1$ and $v_2$ as its endpoints.  The \emph{connected components} of $G$ are the maximal induced subgraphs of $G$ which are connected.  The number of connected components of $G$ will be denoted by $c(G)$.
  
  The \emph{complete graph} $K_n$ on $n$ vertices is the unique simple graph having all possible edges, that is, $E(K_n) = \binom{V}{2}$ where $V = V(K_n)$. The \emph{path graph} $P_n$ is the graph that consists of only an $n$-vertex path, and the \emph{cycle graph} $C_n$ is the graph that consists only of an $n$-vertex cycle.

  Given a graph $G$, there are two commonly used operations that produce new graphs. One is \emph{deletion}: given an edge $e \in E(G)$, the graph of $G$ \emph{with} $e$ \emph{deleted} is the graph $G' = (V(G), E(G) \bk \{e\})$, and is denoted $G \bk e$ or $G-e$. Likewise, if $S$ is a multiset of edges, we use $G \bk S$ or $G-S$ to denote the graph $(V(G),E(G) \bk S)$.

  The other operation is the \emph{contraction of an edge} $e = v_1v_2$, denoted $G / e$.  If $v_1 = v_2$ ($e$ is a loop), we define $G / e = G \bk e$.  Otherwise, we create a new vertex $v^*$, and define $G / e$ as the graph $G'$ with $V(G') = (V(G) \bk \{v_1,v_2\}) \cup v^*$, and $E(G') = (E(G) \bk E(v_1, v_2)) \cup E(v^*)$, where $E(v_1,v_2)$ is the set of edges with at least one of $v_1$ or $v_2$ as an endpoint, and $E(v^*)$ consists of each edge in $E(v_1,v_2) \bk e$ with the endpoint $v_1$ and/or $v_2$ replaced with the new vertex $v^*$.  Note that this is an operation on a (possibly nonsimple) graph that identifies two vertices while keeping and/or creating multi-edges and loops.

  Let $G = (V(G),E(G))$ be a (not necessarily simple) graph. A map $\kappa: V(G) \rightarrow \mathbb{N}_{>0}$ is called a \emph{coloring} of $G$. This coloring is called \emph{proper} if $\kappa(v_1) \neq \kappa(v_2)$ for all $v_1,v_2$ such that there exists an edge $e = v_1v_2$ in $E(G)$. The \emph{chromatic symmetric function} $X_G$ of $G$ is defined as \cite{stanley}
  
  $$X_G(x_1,x_2,\dots) = \sum_{\kappa \text{ proper}} \prod_{v \in V(G)} x_{\kappa(v)} = \sum_{\pi \text{ stable}} \tm_{\lm(\pi)}$$ 
  where the first sum ranges over all proper colorings $\kappa$ of $G$, and the second sum ranges over all (set) partitions of $V(G)$ into stable sets. Note that if $G$ contains a loop then $X_G = 0$, and $X_G$ is unchanged by replacing each multi-edge by a single edge. 
  
  The \emph{Tutte (or bad-coloring) symmetric function} $XB_G$ is defined as an element of $\Lambda[t]$ (the ring of symmetric functions with coefficients in $\mathbb{C}[t]$) as \cite{tutteCrew, stanley} 
  $$
  XB_G(t;x_1,x_2,\dots) = \sum_{\kappa} (1+t)^{e(\kappa)}\prod_{v \in V(G)} x_{\kappa(v)} = \sum_{\pi} (1+t)^{e(\pi)}\tm_{\lm(\pi)}
  $$
  where the first sum ranges over \emph{all} colorings of $G$ (not just the proper ones) letting $e(\kappa)$ be the number of edges of $G$ whose endpoints receive the same color from $\kappa$, and the second sum ranges over \emph{all} (set) partitions of $V(G)$ (not just the stable ones), letting $e(\pi)$ be the number of edges of $G$ whose endpoints lie in the same block of $\pi$.
  
  Note that unlike $X_G$, the Tutte symmetric function is affected by multi-edges and is not annihilated by loops. Furthermore, it is easy to verify that setting $t = -1$ in the Tutte symmetric function recovers the chromatic symmetric function.
  
\subsection{Vertex-Weighted Graphs and their Colorings}

A \emph{vertex-weighted graph} $(G,w)$ consists of a graph $G$ and a weight function $w: V(G) \rightarrow \mathbb{N}_{>0}$. For any $S \subseteq V(G)$ we will denote $w(S) = \sum_{v \in S} w(v)$. 


Given a vertex-weighted graph $(G,w)$ and a non-loop edge $e = v_1v_2 \in E(G)$ we define its \emph{contraction by e} to be the graph $(G/e,w/e)$, where $w/e$ is the weight function such that $(w/e)(v) = w(v)$ if $v$ is not the contracted vertex $v^*$, and $(w/e)(v^*) = w(v_1) + w(v_2)$ (if $e$ is a loop, we define the contraction of $(G,w)$ by $e$ to be $(G \bk e, w)$).

The chromatic symmetric function may be extended to vertex-weighted graphs as
$$
X_{(G,w)} = \sum_{\kappa \text{ proper}} \prod_{v \in V(G)} x_{\kappa(v)}^{w(v)} = \sum_{\pi \text{ stable}} \tm_{\lm(\pi)}
$$
where again the sum ranges over all proper colorings $\kappa$ of $G$, and $\lm(\pi)$ is the integer partition whose parts are $\{w(\pi_i): \pi_i \text{ is a block of } \pi\}$. In this setting the chromatic symmetric function admits the deletion-contraction relation\footnote{This deletion-contraction relation was used in equivalent form for the Hopf algebra of vertex-weighted graphs by Chmutov, Duzhin, and Lando \cite{chmutov}, and for the $W$-polynomial by Noble and Welsh \cite{noble}.} \cite{delcon}
\begin{equation}\label{eq:delcon}
X_{(G,w)} = X_{(G \bk e, w)} - X_{(G/e,w/e)}.
\end{equation}

Similarly, the Tutte symmetric function may be extended as 
$$
XB_{(G,w)} = \sum_{\kappa} (1+t)^{e(\kappa)}\prod_{v \in V(G)} x_{\kappa(v)}^{w(v)} = \sum_{\pi} (1+t)^{e(\pi)}\tm_{\lm(\pi)}
$$
and it may be shown that this function satisfies the deletion-contraction relation \cite{tutteCrew}
\begin{equation}\label{eq:tuttedelcon}
XB_{(G,w)} = XB_{(G \bk e, w)} + tXB_{(G/e,w/e)}.
\end{equation}


\section{Characterizing the Modular Relations of the Chromatic and Tutte Symmetric Functions}

Consider $\mathbb{G}$ as the Hopf algebra of graphs with vertex set $[n]$ for some $n \geq 1$ (taken together with a zero graph). Then $\mathbb{G}$ is a vector space over $\mathbb{C}$ freely generated by these graphs, and its elements are formal sums of such labelled graphs. By extending it linearly, we may view the chromatic symmetric function $X$ as a linear transformation from $\mathbb{G}$ to $\Lambda$. Thus, determining the extent to which the chromatic symmetric function distinguishes graphs is essentially equivalent to determining the kernel of the map $X$. 

In \cite{raul}, Penaguiao determined exactly which elements of $\mathbb{G}$ generate the kernel of $X: \mg \rightarrow \Lambda$. In the next section, we provide a similar argument to find a generating set for the kernel of $XB: \mathbb{G}[t] \rightarrow \Lambda[t]$ as a map\footnote{Here $\mathbb{G}[t]$ consists of (equivalently) polynomials in $t$ with coefficients in $\mathbb{G}$, or formal linear combinations of vertex-labelled graphs with coefficients in $\mathbb{C}[t]$.}. 

Before this, we first introduce in this section new terminology and structural theorems that will aid in our description and its proof. For brevity, if $G$ is a graph and $E$ is a set of edges whose endpoints are a subset of $V(G)$, we will use $G \uplus E$ to mean $(V(G), E(G) \uplus E)$, where $\uplus$ represents the disjoint union of multisets in which we add all copies of each element from both sets (so if an edge $uv$ occurs $m$ times in $E(G)$ and $n$ times in $E$, it occurs $m+n$ times in $E(G) \uplus E$).

\begin{definition}
Let $H_1, \dots, H_k$ be graphs with vertex set $[n]$. Given a linear combination $L = c_1(t)H_1 + \dots + c_k(t)H_k \in \mathbb{G}[t]$ and a graph $G$ on $[N]$ with $N \geq n$, we define the \emph{extension of $L$ by $G$} as 
$$
Ext(L;G) = c_1(t)(G \uplus E(H_1)) + \dots + c_k(t)(G \uplus E(H_k)).
$$
\end{definition}

\begin{definition}
A linear combination $L \in \mathbb{G}[t]$ of graphs with vertex set $[n]$ is a \emph{modular relation for $XB$} if for every graph $G$ with vertex set $[N]$ for $N \geq n$, we have that $Ext(L;G) \in Ker(XB)$.
\medskip

Analogously, a linear combination $l \in \mathbb{G}$ of graphs with vertex set $[n]$ is a \emph{modular relation for $X$} if for every graph $G$ with vertex set $[N]$ for $N \geq n$, we have that $Ext(l;G) \in Ker(X)$.
\end{definition}

The intuitive method used by Penaguaio and used in this paper to determine the kernel of a function $f$ from a graph algebra to a vector space is to use local modifications of certain small (constant-size) subgraphs of any larger graph $G$ to give a combination $\sum c_iG_i$ of graphs such that $f(G) = f(\sum c_iG_i)$, regardless of the size or structure of the rest of the graph $G$ (or how it connects to the subgraphs). It then suffices to find enough such modifications to transform any linear combination of graphs to one in which each graph lies in the preimage of some fixed basis of the image space. The concept of modular relations gives a formal framework for describing all such local graph modifications. Thus, to aid in finding a generating set for the kernel of $XB$, we proceed now to give a structural characterization of its modular relations. We will need an auxiliary definition:

\begin{definition}

An element $L \in \mathbb{G}[t]$ is said to be written in \emph{Tutte standard form} when we write it as
$$
L = \sum_{i=1}^k c_i(1+t)^{n_i}H_i
$$
with $c_i \in \mathbb{C}$, $n_i$ a nonnegative integer, and $H_i$ a graph for all $i$ such that for $i_1 \neq i_2$, $(n_{i_1}, H_{i_1}) \neq (n_{i_2}, H_{i_2})$.

\end{definition}

Note that some of the $H_i$ can be identical. It is clear that each $L$ may be written uniquely in Tutte standard form (up to the order of the terms). We now characterize modular relations of $XB$\footnote{The following theorem and proof may be easily modified to accommodate vertex weights, although we will not need this in what follows.}:

\begin{theorem}\label{thm:tfren}

Let $L = \sum_{i=1}^k c_i(1+t)^{n_i}H_i$ be an element of $\mathbb{G}[t]$ with each $H_i$ a graph with vertex set $[n]$, written in Tutte standard form. Then $L$ is a modular relation for $XB$ if and only if for every partition $\pi \vdash [n]$ we have
\begin{equation}\label{eq:tfren}
B(L;\pi) \coloneqq \sum_{i=1}^k c_i(1+t)^{n_i+e_{H_i}(\pi)} = 0
\end{equation}
where we recall that $e_{H_i}(\pi)$ is the number of edges of $H_i$ with both endpoints in the same block of $\pi$.
\end{theorem}

\begin{proof}

First, we show that if $B(L;\pi) = 0$ for all $\pi \vdash [n]$, then $L$ is a modular relation for $XB$. Let $G$ be any graph with vertex set $[N]$ for $N \geq n$. For $\pi \vdash [N]$, let $\pi' \vdash [n]$ be the partition whose blocks are the nonempty intersections of blocks of $\pi$ with $[n]$. Then
$$XB(Ext(L;G)) = \sum_{i=1}^k c_i\sum_{\pi \vdash [N]} (1+t)^{n_i+e_{G \uplus E(H_i)}(\pi)}\tm_{\lm(\pi)} = \sum_{i=1}^k c_i\sum_{\pi \vdash [N]} (1+t)^{n_i+e_G(\pi)+e_{H_i}(\pi')}\tm_{\lm(\pi)}
$$
$$
 = \sum_{\pi \vdash [N]} (1+t)^{e_G(\pi)}\tm_{\lm(\pi)}\left(\sum_{i=1}^k c_i (1+t)^{n_i+e_{H_i}(\pi')}\right) = 0
$$
so $Ext(L;G) \in Ker(XB)$ for all $G$, and so $L$ is a modular relation for $XB$.

Now, we show that if $B(L;\pi) = 0$ does not hold for all $\pi \vdash [n]$, then $L$ is not a modular relation for $XB$. Thus it suffices to find a graph $G$ with $V(G) = [N]$ for some $N \geq n$ such that $Ext(L;G) \notin Ker(XB)$. In what follows, let $\pi^* = \pi_1^* \sqcup \dots \sqcup \pi_{l(\pi^*)}^*$ be a particular choice of $\pi \vdash [n]$ such that $B(L;\pi^*) \neq 0$ and let $a$ be some particular nonnegative integer such that $[(1+t)^a]B(L;\pi^*) \neq 0$. Let $M = l(\pi^*)(1+\max_i \{n_i+e_{H_i}(\pi^*)\})$, and note that $a < \frac{M}{l(\pi^*)}$. We construct the counterexample graph $G$ as follows:

The vertex set of $G$ will be $[n+Ml(\pi^*)]$. We will view $V(G)$ as the disjoint union of the base vertex set $[n]$ and \emph{clouds} $C_1 , \dots, C_{l(\pi^*)}$, where $C_i = [n+M(i-1)+1, n+Mi]$. The edges of $E(G)$ consist of a single instance each of
\begin{itemize}
\item All possible edges connecting a vertex of $C_i$ to a vertex of $\pi_j^*$ with $i \neq j$.
\item All possible edges connecting a vertex of $C_i$ to a vertex of $C_j$ with $i \neq j$.
\end{itemize}

We claim that $Ext(L;G) \notin Ker(XB)$. More specifically, if $\lm^* = (\lm(\pi^*)_1+M, \dots, \lm(\pi^*)_{l(\pi^*)}+M)$, we claim that
$[(1+t)^a\tm_{\lm^*}]XB(Ext(L;G)) \neq 0$. We evaluate this as
$$
[(1+t)^a\tm_{\lm^*}]XB(Ext(L;G)) = \sum_{i=1}^k c_iN_i(a)
$$
where $N_i(a)$ is the number of partitions $\pi \vdash [n+Ml(\pi^*)]$ with $\lm(\pi) = \lm^*$ such that $n_i + e_{G \uplus H_i}(\pi) = a$. 

It suffices to show that for each $i$ the only partition $\pi \vdash [n+Ml(\pi^*)]$ with type $\lm^*$ and $e_{G \uplus H_i}(\pi) \leq a$ is the partition with blocks $\pi_1^* \cup C_1, \dots, \pi_{l(\pi^*)}^* \cup C_{l(\pi^*)}$, since if so
$$
[(1+t)^a\tm_{\lm^*}]XB(Ext(L;G)) = \sum_{i=1}^k c_iN_i(a) = [(1+t)^a]B(L; \pi^*) \neq 0.
$$

Thus, we fix an arbitrary $i$ and show the claim. Let us consider a fixed partition $p \vdash [n+Ml(\pi^*)]$ with type $\lm^*$ and $e_{G \uplus H_i}(p) \leq  a$. The vertices of the cloud $C_1$ are distributed into the $l(\pi^*)$ blocks of $p$, so some block $B_1$ contains at least $M/l(\pi^*) > a$ vertices of $C_1$. If this block $B_1$ contained any vertex $v \in C_i$ with $i \neq 1$, then as $v$ is connected with every vertex in $C_1$, we would have $e_{G \uplus H_i}(p) > a$, a contradiction. Thus, $B_1$ contains no vertices from any cloud other than $C_1$. Likewise, the vertices of the cloud $C_2$ are distributed in such a way that some block $B_2$ of $\pi$ contains at least $M/l(\pi^*)$ vertices of $C_2$. Clearly $B_1 \neq B_2$, and using the same argument as above we find that $B_2$ cannot contain vertices from any cloud other than $C_2$. Continuing in this manner we find that our partition $p$ must consist of block $B_1, \dots, B_{l(\pi^*)}$ such that all vertices of $C_i$ are in $B_i$ for each $i$. Thus, the partition $p$ necessarily filters the clouds each into their own block.

Now, if $v \in \pi_i^*$, and $v \in B_j$ with $i \neq j$, then since $v$ is connected to every vertex of $C_j$, we have $e_{G \uplus H_i}(p) \geq M > a$, a contradiction, so in fact for each $i$, all vertices of $\pi_i^*$ must go into $B_i$. Thus, the only possibility for the partition $p$ is exactly the one with blocks $\pi_1^* \cup C_1, \dots, \pi_{l(\pi^*)}^* \cup C_{l(\pi^*)}$, so we are done.

\end{proof}

We may also show the following:

\begin{cor}\label{cor:xfren}

Let $l = \sum_{i=1}^k c_iH_i$ be an element of $\mathbb{G}$ with each $H_i$ a graph with vertex set $[n]$ for some positive integer $n$. Then $l$ is a modular relation for $X$ if and only if for every partition $\pi \vdash [n]$,
$$
C(l; \pi) = \sum_{i=1}^k c_i\delta_{e_{H_i}(\pi),0} = 0.
$$

\end{cor}

\begin{proof}

If $C(l;\pi) = 0$ for all $\pi$, then for any graph $G$ on $[N]$ with $N \geq n$ we have (again denoting by $\pi' \vdash [n]$ the partition whose blocks are the nonempty intersections of the blocks of $\pi \vdash [N]$ with $[n]$)
$$
X(Ext(l;G)) = \sum_{i=1}^k c_i \sum_{\pi \vdash [N]} (\delta_{e_{G \uplus E(H_i)}(\pi),0})\tm_{\lm(\pi)} = \sum_{i=1}^k c_i \sum_{\substack{\pi \vdash [N] \\ e_G(\pi) = 0}} (\delta_{e_{H_i}(\pi'),0})\tm_{\lm(\pi)}$$

$$ = \sum_{\substack{\pi \vdash [N] \\ e_G(\pi) = 0}}\tm_{\lm(\pi)} \left(\sum_{i=1}^k c_i\delta_{e_{H_i}(\pi),0}\right) = 0
$$
so $l$ is a modular relation for $X$.  

Conversely, if there exists $\pi$ such that $C(l;\pi) \neq 0$, then we may construct $G$ as in the proof of Theorem \ref{thm:tfren} and verify using the same arguments that $X(Ext(l;G)) \neq 0$, so $l$ is not a modular relation for $X$.
\end{proof}

\begin{cor}\label{cor:tuttexfren}

Given $L \in \mathbb{G}[t]$ such that all graphs in $L$ have vertex set $[n]$, let $L^* \in \mathbb{G}$ be obtained by setting $t = -1$ in $L$. Then if $L$ is a modular relation for $XB$, $L^*$ is a modular relation for $X$.

\end{cor}

\begin{proof}

This follows since for any $\pi \vdash [n]$, if $B(L;\pi) = 0$ in the statement of Theorem \ref{thm:tfren}, then also $C(L^*;\pi) = 0$ in the statement of Corollary \ref{cor:xfren}.

\end{proof}

Among other things, these theorems allow us to formally demonstrate that there is no single local relation that preserves the chromatic symmetric function:

\begin{cor}\label{cor:mod}

Let $H_1$ and $H_2$ be distinct simple graphs with vertex set $[n]$. Then $H_1-H_2$ is never a modular relation for $X$ (and thus not for $XB$).

\end{cor}

\begin{proof}

As $H_1$ and $H_2$ are distinct, there exists an edge $ij$ and $k \in \{1,2\}$ such that $ij \in H_k$ and $ij \notin H_{3-k}$. Then $C(H_1-H_2; p(ij)) \neq 0$ (where we recall $p(ij)$ is the partition of $[n]$ whose only nonsingleton block is $\{i,j\}$), so we are done by Corollary \ref{cor:xfren}.

\end{proof}

Thus, the best we can do to get a two-term ``modular relation" for either $X$ or $XB$ is an argument akin to Theorem 4.2 from \cite{ore} or Theorem 11 from \cite{tutteCrew} that gives such relations for graphs with some guaranteed nice structure (in these cases, a symmetry). Note also that Corollary \ref{cor:mod} holds and gives a nonzero element of $\mathbb{G}[t]$ even when $H_1$ and $H_2$ are isomorphic. Since in this case clearly $H_1-H_2 \in Ker(XB)$, this immediately implies that the set of modular relations of $XB$ is a \emph{proper} subset of $Ker(XB)$ (and likewise for $X$).

\section{The Kernel of $XB$}

As an illustration of how to use the new terminology of Section 3 and to highlight the similarities to our later characterization of $Ker(XB)$, we now provide Penaguaio's characterization of $Ker(X)$ \cite{raul}.

\begin{definition}
Let $G$ be a labelled graph.
\begin{itemize}
 \item If $H$ is a labelled graph such that $G$ and $H$ are isomorphic as unlabelled graphs, define
$$
\ell_{iso}(G,H) = G - H.
$$
\item Let $G$ be the graph with $V(G) = [3]$, and $E(G) = \{12,13,23\}$. Define
$$\ell_{os} = G - G \bk \{12\} - G \bk \{13\} +G \bk \{12,13\}.$$
\begin{figure}[hbt]
\begin{center}
  \begin{tikzpicture}[scale=1.5]
    \node[label=below:{2}, fill=black, circle] at (0, 0)(1){};
    \node[label=below:{3}, fill=black, circle] at (1, 0)(2){};
    \node[label=above:{1}, fill=black, circle] at (0.5, 0.866)(3){};

    \draw[black, thick] (3) -- (1);
    \draw[black, thick] (3) -- (2);
    \draw[black, thick] (2) -- (1);

    \draw (1.5, 0.5) coordinate (MI) node[right] { $\bf{-}$ };

    \node[label=below:{2}, fill=black, circle] at (2, 0)(4){};
    \node[label=below:{3}, fill=black, circle] at (3, 0)(5){};
    \node[label=above:{1}, fill=black, circle] at (2.5, 0.866)(6){};
    
    \draw[black, thick] (5) -- (4);
    \draw[black, thick] (6) -- (4);
    
    \draw (3.5, 0.5) coordinate (MI2) node[right] { $\bf{-}$ };

    \node[label=below:{2}, fill=black, circle] at (4, 0)(7){};
    \node[label=below:{3}, fill=black, circle] at (5, 0)(8){};
    \node[label=above:{1}, fill=black, circle] at (4.5, 0.866)(9){};

    \draw[black, thick] (7) -- (8);
    \draw[black, thick] (8) -- (9);
    
    \draw (5.5, 0.5) coordinate (PL) node[right] { $\bf{+}$ };

    \node[label=below:{2}, fill=black, circle] at (6, 0)(10){};
    \node[label=below:{3}, fill=black, circle] at (7, 0)(11){};
    \node[label=above:{1}, fill=black, circle] at (6.5, 0.866)(12){};

    \draw[black, thick] (10) -- (11);

  \end{tikzpicture}
\end{center}
\label{fig:rel2}
\caption{$\ell_{os}$}
\end{figure}
\end{itemize}
Furthermore, let $T_{iso} \subseteq \mathbb{G}$ be the set of all elements that may be written as $\ell_{iso}(G,H)$ for some graphs $G$ and $H$, and let $T_{os} \subseteq \mathbb{G}$ be the set of all extensions $Ext(\ell_{os};G)$ for some $G$.
\end{definition}

The relation $\ell_{os}$ is so named because it was originally discovered by Orellana and Scott as the first known modular relation for the chromatic symmetric function \cite{ore} (which is easy to verify using Corollary \ref{cor:xfren}). Penaguaio shows that any single graph in $\mathbb{G}$ (with coefficient $1$) may be expressed as a linear combination of complete multipartite graphs by repeatedly adding only elements of $T_{iso}$ and $T_{os}$ \cite{raul}. Since the set of chromatic symmetric functions of these graphs forms a basis for $\Lambda$ (see \cite{alif}), this is sufficient to show the following:

\begin{theorem}[\cite{raul}, Theorem 2]

Let $I_X \subseteq \mathbb{G}$ be the set of all linear combinations of elements of $T_{iso}$, and let $OS$ be the set of all linear combinations of elements of $T_{os}$. Then $span(I_X, OS) = Ker(X)$.

\end{theorem}

We give an analogous result for the Tutte symmetric function using the basis of \emph{star forests}, which we define here. 

The \emph{star} $S_k$ on $k$ vertices is the simple connected graph where one vertex has degree $k-1$, and all other vertices have degree $1$. A star with three or more vertices is said to be \emph{rooted} at its vertex of largest degree. If $\lm$ is an integer partition, the \emph{star forest} $S_{\lm}$ is the graph with connected components $S_{\lm_1}, \dots, S_{\lm_{l(\lm)}}$. It is known (e.g. by \cite{alif}) that $\{X(S_{\lm}) | \lm \text{ an integer partition}\}$ is a basis for $\Lambda$.

In our proof, we will need the following auxiliary lemma:

\begin{lemma}\label{lem:star}

Let $G$ be a simple graph with vertex set $[n]$. Suppose that for any three distinct vertices $a < b < c$ of $G$, we have either

\begin{itemize}
    \item $\{ab, ac, bc\} \cap E(G) \leq 1$, or
    \item $\{ab, ac, bc\} \cap E(G) = \{ac, bc\}$.
\end{itemize}

Then $G$ is a star forest, and furthermore, each induced star $S$ of $G$ with three or more vertices is rooted at its largest vertex.

\end{lemma}

\begin{proof}

Note that our condition ensures that no vertex has two distinct neighbors that are greater than it.  Thus, $G$ cannot contain a cycle, since then the least vertex in the cycle would violate this. We may conclude then that $G$ is a forest.

Now, let $C$ be an arbitrary connected component of $G$, and let $v_0$ be the largest vertex in $C$. Suppose that there exists $v \in V(C)-v_0$ such that $vv_0 \notin E(G)$. Then, as $C$ is a tree, there is a unique path from $v$ to $v_0$ containing at least one vertex other than $v$ and $v_0$. Let $v_1$ be the vertex in this path that is adjacent to $v_0$, and let $v_2$ be the other vertex of the $vv_0$ path that $v_1$ is adjacent to (it may be $v$). Then the triple $\{v_1, v_2, v_0\}$ violates our condition, since $\{v_1v_2, v_1v_0, v_2v_0\} \cap E(G)$ contains at least two edges, and they are not the two edges incident to the largest vertex $v_0$. This is a contradiction, and it follows that for every $v \in V(C)-v_0$, we have $vv_0 \in E(G)$. Since $C$ is a tree, this accounts for all edges of $C$, and the result follows.

\end{proof}

Since we will use star forests of the specific type given in Lemma \ref{lem:star} throughout our next proof, we will give them a name.  A graph $G$ with vertex set $[n]$ is called a \emph{bright} star forest if every induced star is rooted at its largest vertex. If $G$ is not a bright star forest, it is called \emph{dull}.

We now introduce the modular relations that will generate $Ker(XB)$ (all defined terms are elements of $\mathbb{G}[t]$):

\begin{definition}
\begin{itemize}
\item Let $G_1$ be the graph with vertex set $[1]$ and edge set $\{11\}$ (so $G_1$ contains just a loop at this vertex). Define $$\ell_{loop} = G_1-(t+1)G_1 \bk \{11\}$$
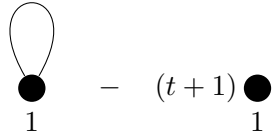
\begin{figure}[hbt]
\begin{center}
  \begin{tikzpicture}[scale=1.5]
    \node[label=below:{1}, fill=black, circle] at (0, 1)(1){};

    \path (1) edge [loop above] node {} (1);

    \draw (0.5, 1) coordinate (MI) node[right] { $\bf{-}$ };

    \draw (1, 1) coordinate (T2) node[right] { $(t+1)$ };

    \node[label=below:{1}, fill=black, circle] at (2, 1)(3){};
    
    \end{tikzpicture}
\end{center}
\label{fig:rel1}
\caption{$\ell_{loop}$}
\end{figure}
\item Let $G_2$ be the graph with vertex set $[2]$ and edge set $\{e_1 = e_2 = 12\}$ (so $G$ has two edges between two vertices). Define 
$$\ell_{multi} = G_2-(t+2)G_2 \bk \{e_2\}+(t+1)G_2 \bk \{e_1, e_2\}.$$
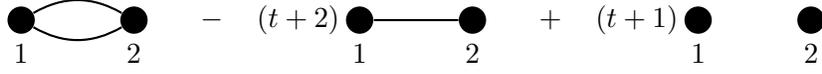
\begin{figure}[hbt]
\begin{center}
  \begin{tikzpicture}[scale=1.5]
    \node[label=below:{1}, fill=black, circle] at (0, 1)(1){};
    \node[label=below:{2}, fill=black, circle] at (1, 1)(2){};

    \draw[black, thick] (1) edge [bend left=30] (2);
    \draw[black, thick] (1) edge [bend right=30] (2);

    \draw (1.5, 1) coordinate (MI) node[right] { $\bf{-}$ };

    \draw (2, 1) coordinate (T2) node[right] { $(t+2)$ };

    \node[label=below:{1}, fill=black, circle] at (3, 1)(3){};
    \node[label=below:{2}, fill=black, circle] at (4, 1)(4){};
    \draw[black, thick] (3) -- (4);

    \draw (4.5, 1) coordinate (MI) node[right] { $\bf{+}$ };
    \draw (5, 1) coordinate (T1) node[right] { $(t+1)$ };

    \node[label=below:{1}, fill=black, circle] at (6, 1)(5){};
    \node[label=below:{2}, fill=black, circle] at (7, 1)(6){};
    
    \end{tikzpicture}
\end{center}
\label{fig:lm}
\caption{$\ell_{multi}$}
\end{figure}

\item Let $G_3$ be the graph with vertex set $[3]$ and edge set $\{12,23\}$. Define
$$\ell_{os+} = G_3 + G_3^c - G_3 \bk \{23\} - (G_3 \bk \{23\})^c.$$
\begin{figure}[hbt]
\begin{center}
  \begin{tikzpicture}[scale=1.5]
    \node[label=below:{2}, fill=black, circle] at (0, 0)(1){};
    \node[label=below:{3}, fill=black, circle] at (1, 0)(2){};
    \node[label=above:{1}, fill=black, circle] at (0.5, 0.866)(3){};

    \draw[black, thick] (3) -- (1);
    \draw[black, thick] (2) -- (1);

    \draw (1.5, 0.5) coordinate (PL) node[right] { $\bf{+}$ };
    
    \node[label=below:{2}, fill=black, circle] at (2.5, 0)(4){};
    \node[label=below:{3}, fill=black, circle] at (3.5, 0)(5){};
    \node[label=above:{1}, fill=black, circle] at (3, 0.866)(6){};
    
    \draw[black, thick] (5) -- (6);
    
    \draw (4, 0.5) coordinate (MI) node[right] { $\bf{-}$ };
    
    \node[label=below:{2}, fill=black, circle] at (5, 0)(7){};
    \node[label=below:{3}, fill=black, circle] at (6, 0)(8){};
    \node[label=above:{1}, fill=black, circle] at (5.5, 0.866)(9){};
    
    \draw[black, thick] (8) -- (7);
    \draw[black, thick] (9) -- (8);
    
    \draw (6.5, 0.5) coordinate (M2) node[right] { $\bf{-}$ };
    
    \node[label=below:{2}, fill=black, circle] at (7.5, 0)(A){};
    \node[label=below:{3}, fill=black, circle] at (8.5, 0)(B){};
    \node[label=above:{1}, fill=black, circle] at (8, 0.866)(C){};
    
    \draw[black, thick] (A) -- (C);
  
  \end{tikzpicture}
\end{center}
\label{fig:ls}
\caption{$\ell_{os+}$}
\end{figure}
\end{itemize}
Furthermore, let $T_{loop} \subseteq \mathbb{G}[t]$ be the set of all elements that are equal to $Ext(\ell_{loop};G)$ for some choice of $G$ and likewise for $T_{multi}$ and $T_{os+}$.
\end{definition}

Note that unlike in the case of simple graphs, since we are considering these relations for multigraphs, we only require the edges listed to be present, but there can be multiple copies of these edges, as well as possibly other edges not among those listed. In particular, it is worth noting that the relation given by $\ell_{os+}$ is quite powerful, since it may be applied whenever there are three vertices with \emph{at least} two (distinct) edges among them; all three distinct edges may be present if we don't restrict to simple graphs. In fact, the extension of $\ell_{os+}$ to the complete graph on the vertex set $[3]$ with all multi-edges reduced to single edges is precisely (up to isomorphism) $\ell_{os}$, so we may view $\ell_{os+}$ as a generalization of this relation that is valid for $XB$, hence the choice of notation.

\begin{theorem}\label{thm:kerxb}
Define subsets $I, S_1, S_2 \subseteq \mathbb{G}[t]$ as follows:
\begin{itemize}
\item Let $I$ be the set of all linear combinations of elements of $T_{iso}$. 
\item Let $S_1$ be the set of all linear combinations of elements of $T_{loop} \cup T_{multi}$.
\item Let $S_2$ be the set of all linear combinations of elements of $T_{os+}$.
  \end{itemize}
Then $Span(I, S_1, S_2) = Ker(XB)$.
\end{theorem}

\begin{proof}

First, it is easy to verify that $Span(I, S_1, S_2) \subseteq Ker(XB)$. Clearly by definition $I \subseteq Ker(XB)$, and using Theorem \ref{thm:tfren} every element of $S_1$ and $S_2$ may be seen to be a modular relation for $XB$.

It remains to show that $Ker(XB) \subseteq Span(I, S_1, S_2)$. We will do so by taking an arbitrary element $k \in Ker(XB)$ and show that there exist $i \in I, s_1 \in S_1, s_2 \in S_2$ such that $k - i - s_1 - s_2 = 0$, from which it will follow that $k = i+s_1+s_2 \in Span(I, S_1, S_2)$.

The key observation is that for every single graph $G \in \mathbb{G}[t]$ (with coefficient $1$), we may add elements of $I$, $S_1$, and $S_2$ to it to turn it into a linear combination of bright star forests. We will do this in a step-by-step process. First, if $G$ has a loop at a vertex $v$, we add an element of $I$ to relabel and get a graph $G_0$ with a loop at $1$, and then subtract $Ext(\ell_{loop};G_0)$ to obtain a linear combination of graphs with no loop at $1$ (from here on out, we will simply use ``relabel $G$" to formally mean subtracting an appropriate element of $I$ as here, and simply keeping the same name $G$ for the resulting graph for brevity). If $G$ has other loops, continue this process of subtracting terms from $T_{iso}$ and $T_{loop}$ until $G$ is expressed as a linear combination $L_1$ of loopless graphs. Now, for each term $c_1G_1$ in $L_1$ that has multiple edges between two vertices $v_1$ and $v_2$, we relabel $G_1$ and subtract $Ext(\ell_{multi};G_1)$ to obtain a linear combination of loopless graphs with fewer edges between $v_1$ and $v_2$. Continuing in this manner, by only subtracting terms from $I$ and $S_1$, we may reduce $G$ to a linear combination $L_2$ of simple graphs.

Now, we will subtract elements of $I$, $S_1$, and $S_2$ to write $L_2$ as a linear combination of bright star forests. It suffices to show that we can do this for any simple graph (with coefficient $1$).
We will show this by considering a partial order $\preceq$ on the (finite) set of all simple graphs with vertex set $[n]$. For an edge $e$ in such a graph, define its \emph{right endpoint} to be the larger of its two endpoints. Let $r_1(G) ,\dots, r_{|E(G)|}(G)$ be the list of right endpoints of $G$ with multiplicity, where $r_1 \leq \dots \leq r_{|E(G)|}$. Then for two graphs $G_1$ and $G_2$ with vertex set $[n]$, we define $G_1 \preceq G_2$ if and only if 
\begin{itemize}
    \item $G_1$ has more edges than $G_2$, or
    \item $G_1$ and $G_2$ have the same number $E$ of edges, and there exists $i \in [E]$ such that $r_j(G_1) = r_j(G_2)$ for all $j < i$, and $r_i(G_1) < r_i(G_2)$.
\end{itemize}

We claim that any simple graph $H$ with vertex set $[n]$ that is not a bright star forest may be written as a linear combination of graphs $H_1, \dots, H_k$ such that $H \prec H_i$ for all $i \in [k]$. If $H$ is such a graph, then it necessarily has three vertices $h_1 < h_2 < h_3$ that are in violation of the conditions of Lemma \ref{lem:star}. Relabelling without changing the order, we may temporarily let these be $1 < 2 < 3$.

Case 1: $\{12, 13, 23\} \cap E(H) = \{12, 13\}$. Then we subtract the element obtained by relabelling each graph of $Ext(\ell_{os+};H)$ by swapping $1$ and $2$.

\begin{figure}[hbt]
\begin{center}
  \begin{tikzpicture}[scale=1.5]
    \node[label=below:{1}, fill=black, circle] at (0, 0)(1){};
    \node[label=below:{3}, fill=black, circle] at (1, 0)(2){};
    \node[label=above:{2}, fill=black, circle] at (0.5, 0.866)(3){};

    \draw[black, thick] (3) -- (1);
    \draw[black, thick] (2) -- (1);

    \draw (1.5, 0.5) coordinate (PL) node[right] { $\bf{+}$ };
    
    \node[label=below:{1}, fill=black, circle] at (2.5, 0)(4){};
    \node[label=below:{3}, fill=black, circle] at (3.5, 0)(5){};
    \node[label=above:{2}, fill=black, circle] at (3, 0.866)(6){};
    
    \draw[black, thick] (5) -- (6);
    
    \draw (4, 0.5) coordinate (MI) node[right] { $\bf{-}$ };
    
    \node[label=below:{1}, fill=black, circle] at (5, 0)(7){};
    \node[label=below:{3}, fill=black, circle] at (6, 0)(8){};
    \node[label=above:{2}, fill=black, circle] at (5.5, 0.866)(9){};
    
    \draw[black, thick] (8) -- (7);
    \draw[black, thick] (9) -- (8);
    
    \draw (6.5, 0.5) coordinate (M2) node[right] { $\bf{-}$ };
    
    \node[label=below:{1}, fill=black, circle] at (7.5, 0)(A){};
    \node[label=below:{3}, fill=black, circle] at (8.5, 0)(B){};
    \node[label=above:{2}, fill=black, circle] at (8, 0.866)(C){};
    
    \draw[black, thick] (A) -- (C);
  
  \end{tikzpicture}
\end{center}
\label{fig:case1}
\caption{Case 1 - Subtracting the extension of this term has the effect of locally changing the configuration shown in the left graph to a linear combination of the other three.}
\end{figure}
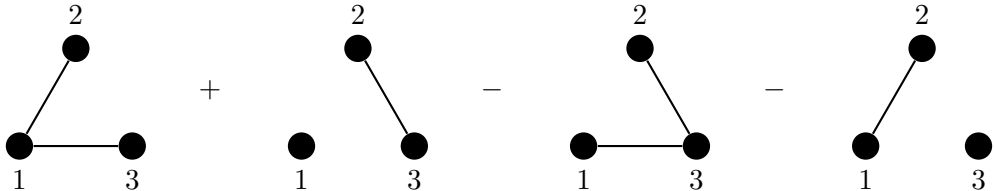

Case 2: $\{12, 13, 23\} \cap E(H) = \{12, 23\}$. Then we subtract $Ext(\ell_{os+};H)$.

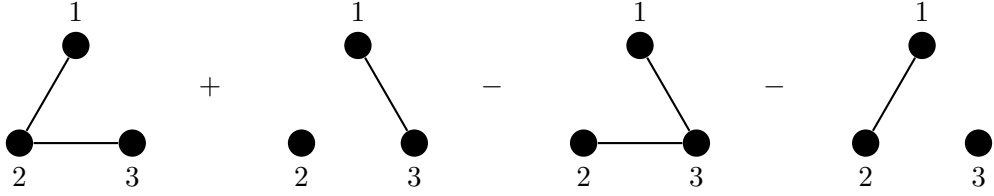
\begin{figure}[hbt]
\begin{center}
  \begin{tikzpicture}[scale=1.5]
    \node[label=below:{2}, fill=black, circle] at (0, 0)(1){};
    \node[label=below:{3}, fill=black, circle] at (1, 0)(2){};
    \node[label=above:{1}, fill=black, circle] at (0.5, 0.866)(3){};

    \draw[black, thick] (3) -- (1);
    \draw[black, thick] (2) -- (1);

    \draw (1.5, 0.5) coordinate (PL) node[right] { $\bf{+}$ };
    
    \node[label=below:{2}, fill=black, circle] at (2.5, 0)(4){};
    \node[label=below:{3}, fill=black, circle] at (3.5, 0)(5){};
    \node[label=above:{1}, fill=black, circle] at (3, 0.866)(6){};
    
    \draw[black, thick] (5) -- (6);
    
    \draw (4, 0.5) coordinate (MI) node[right] { $\bf{-}$ };
    
    \node[label=below:{2}, fill=black, circle] at (5, 0)(7){};
    \node[label=below:{3}, fill=black, circle] at (6, 0)(8){};
    \node[label=above:{1}, fill=black, circle] at (5.5, 0.866)(9){};
    
    \draw[black, thick] (8) -- (7);
    \draw[black, thick] (9) -- (8);
    
    \draw (6.5, 0.5) coordinate (M2) node[right] { $\bf{-}$ };
    
    \node[label=below:{2}, fill=black, circle] at (7.5, 0)(A){};
    \node[label=below:{3}, fill=black, circle] at (8.5, 0)(B){};
    \node[label=above:{1}, fill=black, circle] at (8, 0.866)(C){};
    
    \draw[black, thick] (A) -- (C);
  
  \end{tikzpicture}
\end{center}
\label{fig:case2}
\caption{Case 2 - As above, the configuration in the left-hand graph is turned into a linear combination of those in the other graphs.}
\end{figure}

Case 3: $\{12, 13, 23\} \cap E(H) = \{12, 13, 23\}$. Then we subtract the element obtained by relabelling each graph of $Ext(\ell_{os+};H)$ by swapping $1$ and $2$. We further subtract elements of $I$ and $S_1$ to turn multi-edges into a single edge as applicable.

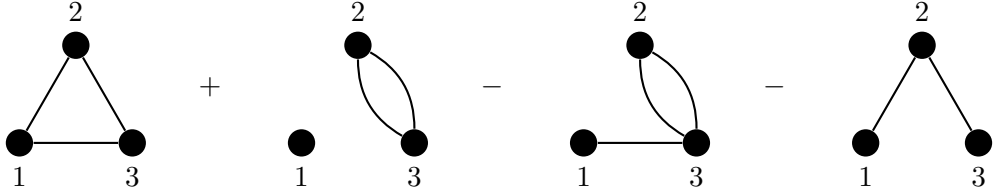
\begin{figure}[hbt]
\begin{center}
  \begin{tikzpicture}[scale=1.5]
    \node[label=below:{1}, fill=black, circle] at (0, 0)(1){};
    \node[label=below:{3}, fill=black, circle] at (1, 0)(2){};
    \node[label=above:{2}, fill=black, circle] at (0.5, 0.866)(3){};

    \draw[black, thick] (3) -- (1);
    \draw[black, thick] (2) -- (1);
    \draw[black, thick] (3) -- (2);

    \draw (1.5, 0.5) coordinate (PL) node[right] { $\bf{+}$ };
    
    \node[label=below:{1}, fill=black, circle] at (2.5, 0)(4){};
    \node[label=below:{3}, fill=black, circle] at (3.5, 0)(5){};
    \node[label=above:{2}, fill=black, circle] at (3, 0.866)(6){};
    
    \draw[black, thick] (5) edge [bend left=30] (6);
    \draw[black, thick] (5) edge [bend right=30] (6);
    
    \draw (4, 0.5) coordinate (MI) node[right] { $\bf{-}$ };
    
    \node[label=below:{1}, fill=black, circle] at (5, 0)(7){};
    \node[label=below:{3}, fill=black, circle] at (6, 0)(8){};
    \node[label=above:{2}, fill=black, circle] at (5.5, 0.866)(9){};
    
    \draw[black, thick] (8) -- (7);
    \draw[black, thick] (9) edge [bend left=30] (8);
    \draw[black, thick] (9) edge [bend right=30] (8);
    
    \draw (6.5, 0.5) coordinate (M2) node[right] { $\bf{-}$ };
    
    \node[label=below:{1}, fill=black, circle] at (7.5, 0)(A){};
    \node[label=below:{3}, fill=black, circle] at (8.5, 0)(B){};
    \node[label=above:{2}, fill=black, circle] at (8, 0.866)(C){};
    
    \draw[black, thick] (A) -- (C);
    \draw[black, thick] (B) -- (C);
  
  \end{tikzpicture}
\end{center}
\label{fig:case3}
\caption{Case 3 - Same as the other two cases, except that we then remove the multi-edges as well.}
\end{figure}

Regardless of which case was used, we relabel again to restore the vertices to $h_1 < h_2 < h_3$. In all cases, we will have subtracted elements of $I$, $S_1$, and $S_2$ from $H$ to obtain a linear combination of three other graphs $H_1$, $H_2$, and $H_3$, and it is easy to check that indeed $H \prec H_1$, $H \prec H_2$ and $H \prec H_3$ all hold. Now, as long as any of the $H_i$ still remain dull, we may apply this process again to each such graph. Since each application of this process produces graphs that are strictly larger in a finite poset, this process must terminate after a finite number of applications, necessarily in a linear combination of bright star forests.

Returning to the main proof, we have now demonstrated that for any $G \in \mathbb{G}[t]$ with vertex set $[n]$, we may subtract elements of $I$, $S_1$, and $S_2$ to write $G$ as a linear combination $L_3$ of bright star forests.

Now, for each integer partition $\lm \vdash n$, we choose a canonical bright star forest with vertex set $[n]$ by defining $R_{\lm}$ to be the star forest such that the stars are induced by the vertex sets $\{1, \dots, \lm_1\}, \{\lm_1+1 \dots, \lm_1+\lm_2\}, \dots, \{|V(G)|-\lm_l+1, \dots, |V(G)|\}$, and each star is rooted at its largest vertex. By subtracting further elements of $I$, we may rewrite $L_3$ so that every graph that occurs is one of the $R_{\lm}$, and thus $L_3$ may be written as  $L_4 = \sum_{\lm, k} c_{\lm, k}(1+t)^kR_{\lm}$ for some $c_{\lm, k} \in \mathbb{C}[t]$.

Now, let $k \in Ker(XB)$ be arbitrary. Then we have shown that there exist $i \in I, s_1 \in S_1, s_2 \in S_2$ such that $k - i - s_1 - s_2$ is a linear combination of terms of the form $(1+t)^kR_{\lm}$, and also $k-i-s_1-s_2 \in Ker(XB)$. By e.g. \cite{alif}, the functions $X(R_{\lm})$ form a basis for $\Lambda$; it is then easy to check that the functions $\{(1+t)^kX(R_{\lm}) | k \in \mathbb{N}, \lm \text{ an integer partition}\}$ form a basis for $\Lambda[t]$. Thus, if $L$ is a finite linear combination of terms of the form $(1+t)^k R_{\lm}$ such that $L \in Ker(XB)$, then $L = 0$. Thus, it must be the case that $k-i-s_1-s_2 = 0$, and this finishes the proof.\footnote{It is possible to use a similar argument to show that considering $X$ and $XB$ as functions from the algebra of labelled, vertex-weighted graphs, their kernels are generated by isomorphisms and those relations arising from the deletion-contraction relations \eqref{eq:delcon} and \eqref{eq:tuttedelcon} respectively \cite{priv}.}
\end{proof}

Essentially, this shows that $Ker(XB)$, after adding generators to remove loops and multi-edges, is generated by isomorphisms and a natural extension of the Orellana-Scott triangular relation to $XB$. 

The modular relation $\ell_{os}$ for $X$ may be generalized to $n$-cycles, as by Dahlberg and van Willigenburg \cite{dahl}, so it is natural to ask if there is a cycle analogue of $\ell_{os+}$ and if it may be generalized. We go a step further, and work with a more general graph class than cycles.

A graph $G$ is called \emph{two-edge-connected} if it is connected, and also $G - e$ is connected for every edge $e \in E(G)$.  Equivalently, whenever there is a partition of $V(G)$ into two nonempty blocks $B_1$ and $B_2$, there are at least two edges that have one endpoint in $B_1$ and one endpoint in $B_2$.

\begin{theorem}\label{thm:twocon}

Let $G$ be a two-edge-connected graph with vertex set $[n]$, and let $e_i$ and $e_j$ be two (not necessarily distinct) edges of $G$. For $S \subseteq E(G)$, let $S$ also be used in a slight abuse of notation to mean the graph with vertex set $[n]$ and edge set $S$. Then 
$$
\sum_{\substack{S \subseteq E(G) \\ e_i \in S}} (-1)^{|S|}S + (1+t)\sum_{\substack{S \subseteq E(G) \\ e_j \notin S }} (-1)^{|S|}S
$$
is a modular relation for $XB$.
\end{theorem}

\begin{proof}

By applying Theorem \ref{thm:tfren}, it suffices to show that for each $\pi \vdash V(G)$, we have that 
$$
\sum_{\substack{S \subseteq E(G) \\ e_i \in S}} (-1)^{|S|}(1+t)^{e_S(\pi)} + \sum_{\substack{S \subseteq E(G) \\ e_j \notin S }} (-1)^{|S|}(1+t)^{e_S(\pi)+1} = 0.
$$

First, we consider the case where $\pi$ has more than one part. Give some arbitrary total ordering to the edges of $G$. By assumption, since $\pi$ has more than one part, there are at least two edges whose endpoints do not lie in the same block of $\pi$, and thus at least one such edge that is not $e_i$. Among these, choose the smallest edge with respect to the ordering; let this edge be $\epsilon$. Let $E_i$ denote the set of subsets of $E(G)$ that contain $e_i$. Then we have an involution $\iota: E_i \rightarrow E_i$ given by $\iota(S) = S \cup \epsilon$ if $S$ does not contain $\epsilon$, and otherwise $\iota(S) = S - \epsilon$. Clearly $|S|$ and $|\iota(S)|$ differ by one, and $e_{S}(\pi) = e_{\iota(S)}(\pi)$, so the first sum may be rewritten as
$$
\sum_{\substack{S \subseteq E(G) \\ e_i, \epsilon \in S}} \left((-1)^{|S|}(1+t)^{e_S(\pi)}+(-1)^{|\iota(S)|}(1+t)^{e_{\iota(S)}(\pi)}\right) $$ $$ = \sum_{\substack{S \subseteq E(G) \\ e_i, \epsilon \in S}} \left((-1)^{|S|}(1+t)^{e_S(\pi)}+(-1)^{|S|- 1}(1+t)^{e_S(\pi)}\right) = 0.
$$

An analogous argument shows that the second sum also goes to $0$ for any $\pi$ that is not equal to the whole vertex set.

It remains to check the case where $\pi$ consists of a single block. In this case the overall sum reduces to
$$
\sum_{\substack{S \subseteq E(G) \\ e_i \in S}} (-1)^{|S|}(1+t)^{|S|} + \sum_{\substack{S \subseteq E(G) \\ e_j \notin S }} (-1)^{|S|}(1+t)^{|S|+1} $$ $$ = \sum_{\substack{S \subseteq E(G) \\ e_i \in S}} (-1)^{|S|}(1+t)^{|S|} + \sum_{\substack{S \subseteq E(G) \\ e_j \in S }} (-1)^{|S|-1}(1+t)^{|S|} = 0.
$$

\end{proof}

By Corollary \ref{cor:tuttexfren}, we also get the following:
\begin{cor}
Let $G$ be a two-edge-connected graph with vertex set $[n]$, and let $e_i$ be any edge of $G$. Then 
$$
\sum_{\substack{S \subseteq E(G) \\ e_i \in S}} (-1)^{|S|}S
$$
is a modular relation for $X$.

\end{cor}

Since cycles are two-edge-connected, these results generalize the $n$-cycle modular relation for the chromatic symmetric function of Dahlberg and van Willigenburg \cite{dahl} (and the proof we give is essentially a modification of their edge-swapping argument). Note that from Theorem \ref{thm:twocon} we also obtain an analogous $n$-cycle relation for the Tutte symmetric function:

\begin{cor}\label{thm:cycle}
Let $n \geq 3$ be a positive integer and let $C$ be an $n$-vertex cycle of a graph $G$. Let $v_1, v_2, \dots, v_n$ be the consecutive vertices of the cycle, and let the edges be $e_1 = v_nv_1, e_2 = v_1v_2, \dots, e_n = v_nv_1$. Then for any $1 \leq i,j \leq n$, we have
\begin{equation}\label{eq:cycle}
\sum_{\substack{S \subseteq E(C) \\ e_i \notin S}} (-1)^{|S|}XB_{G \backslash S} + (1+t)\sum_{\substack{S \subseteq E(C) \\ e_j \in S}} (-1)^{|S|}XB_{G \backslash S} = 0.
\end{equation}
\end{cor}

\begin{figure}[hbt]
\begin{center}
  \begin{tikzpicture}[scale=1.5]
    \node[label=below:{2}, fill=black, circle] at (0, 0)(1){};
    \node[label=below:{3}, fill=black, circle] at (1, 0)(2){};
    \node[label=above:{1}, fill=black, circle] at (0.5, 0.866)(3){};

    \draw[black, thick] (3) -- (1);
    \draw[black, thick] (3) -- (2);
    \draw[black, thick] (2) -- (1);

    \draw (1.5, 0.5) coordinate (MI) node[right] { $\bf{-}$ };

    \node[label=below:{2}, fill=black, circle] at (2, 0)(4){};
    \node[label=below:{3}, fill=black, circle] at (3, 0)(5){};
    \node[label=above:{1}, fill=black, circle] at (2.5, 0.866)(6){};
    
    \draw[black, thick] (5) -- (4);
    \draw[black, thick] (6) -- (4);
    
    \draw (3.5, 0.5) coordinate (MI2) node[right] { $\bf{-}$ };
    \draw (4, 0.5) coordinate (T2) node[right] { $(t+2)$ };

    \node[label=below:{2}, fill=black, circle] at (5, 0)(7){};
    \node[label=below:{3}, fill=black, circle] at (6, 0)(8){};
    \node[label=above:{1}, fill=black, circle] at (5.5, 0.866)(9){};

    \draw[black, thick] (7) -- (8);
    \draw[black, thick] (8) -- (9);
    
  \end{tikzpicture}

  \vspace{0.5cm}

  \begin{tikzpicture}[scale=1.5]
    \draw (0.5, 0.5) coordinate (PL) node[right] { $\bf{+}$ };
    \draw (1, 0.5) coordinate (T2) node[right] { $(t+2)$ };

    \node[label=below:{2}, fill=black, circle] at (2, 0)(1){};
    \node[label=below:{3}, fill=black, circle] at (3, 0)(2){};
    \node[label=above:{1}, fill=black, circle] at (2.5, 0.866)(3){};

    \draw[black, thick] (1) -- (2);

    \draw (3.5, 0.5) coordinate (P2) node[right] { $\bf{+}$ };
    \draw (4, 0.5) coordinate (T1) node[right] { $(t+1)$ };

    \node[label=below:{2}, fill=black, circle] at (5, 0)(4){};
    \node[label=below:{3}, fill=black, circle] at (6, 0)(5){};
    \node[label=above:{1}, fill=black, circle] at (5.5, 0.866)(6){};

    \draw[black, thick] (5) -- (6);

    \draw (6.5, 0.5) coordinate (M3) node[right] { $\bf{-}$ };
    \draw (7, 0.5) coordinate (TE) node[right] { $(t+1)$ };

    \node[label=below:{2}, fill=black, circle] at (8, 0)(7){};
    \node[label=below:{3}, fill=black, circle] at (9, 0)(8){};
    \node[label=above:{1}, fill=black, circle] at (8.5, 0.866)(9){};

  \end{tikzpicture}
\end{center}
\label{fig:nis3}
\caption{The case $n = 3$ of Corollary \ref{thm:cycle}}
\end{figure}
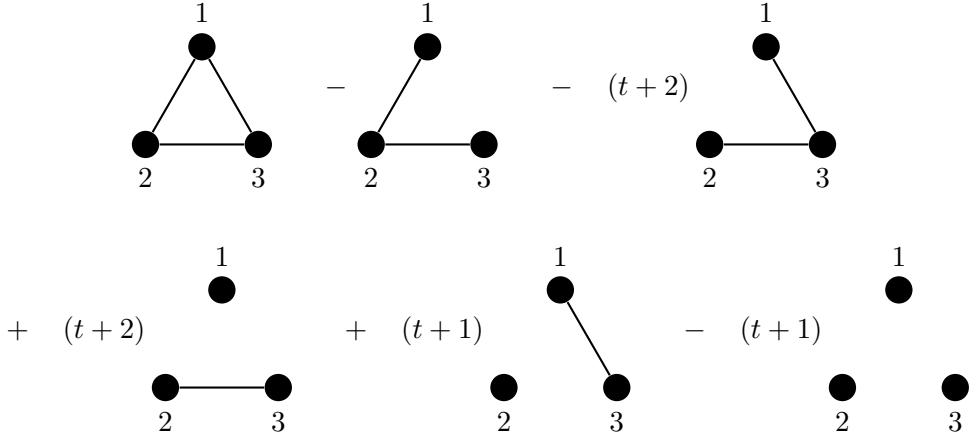

For example, Figure \ref{fig:nis3} illustrates the case $n = 3$. Here $e_i = 23$ and the corresponding sum contributes to the first four terms, and $e_j = 12$ and the corresponding sum contributes to the last four terms.

Note that setting $t = -1$ in \eqref{eq:cycle} also recovers the Dahlberg-van Willigenburg $n$-cycle result \cite{dahl}.

\section{Acknowledgments}

The authors would like to acknowledge Jos\'e Aliste-Prieto and Jos\'e Zamora for introducing us to Penaguaio's work on the kernel of $X$ and for helpful discussions and comments on the kernel of $X$ and $XB$. We would also like to thank the referees for their helpful comments.

We acknowledge the support of the Natural Sciences and Engineering Research Council of Canada (NSERC), [funding
reference number RGPIN-2020-03912]. Cette recherche a \'et\'e financ\'ee par le Conseil de recherches en sciences naturelles et en g\'enie du Canada (CRSNG),
[num\'ero de r\'ef\'erence RGPIN-2020-03912].

Declaration of Interests: none.
\bibliographystyle{plain}
\bibliography{bib}

\end{document}